\input amstex
\input amsppt.sty
\magnification=\magstep1
\hsize=30truecc
\vsize=22.2truecm
\baselineskip=16truept
\TagsOnRight
\nologo
\pageno=1
\topmatter
\def\N{\Bbb N}
\def\Z{\Bbb Z}

\def\l{\left}
\def\r{\right}
\def\b{\bigg}

\def\({\b(}
\def\[{\b[}
\def\){\b)}
\def\]{\b]}

\def\t{\text}
\def\f{\frac}
\def\mo{\roman{mod}}
\def\ord{\roman{ord}}

\def\bi{\binom}
\def\eq{\equiv}

\def\ls{\leqslant}
\def\gs{\geqslant}
\def\al{\alpha}

\def\da{\delta}

\def\Proof{\noindent{\it Proof}}
\def\Remark{\medskip\noindent{\it Remark}}

\def\bm#1#2#3{\thickfracwithdelims ()\thickness0{#1}{#2}_{#3}}

\hbox {Acta Arith. 126(2007), no.\,4, 387--398.}
\bigskip
\title  Combinatorial congruences and Stirling numbers\endtitle
\author {Zhi-Wei Sun (Nanjing)}\endauthor
\leftheadtext{Z. W. Sun}
 \abstract In this paper we obtain some
sophisticated combinatorial congruences involving binomial
coefficients and confirm two conjectures of the author and Davis. They are closely
related to our investigation of the periodicity of the sequence
$\sum_{j=0}^l\bi ljS(j,m)a^{l-j}\ (l=m,m+1,\ldots)$ modulo a prime
$p$, where $a$ and $m>0$ are integers, and those $S(j,m)$ are
Stirling numbers of the second kind. We also give a new extension
of Glaisher's congruence by showing that
$(p-1)p^{\lfloor\log_pm\rfloor}$ is a period of the sequence
$\sum_{j\eq r\,(\mo\ p-1)}\bi ljS(j,m)\ (l=m,m+1,\ldots)$ modulo
$p$.
\endabstract
\thanks 2000 {\it Mathematics Subject Classification}.\,Primary 11B65;
Secondary 05A10, 11A07, 11B73.\newline\indent
Supported by the National Science Fund
for Distinguished Young Scholars (Grant No. 10425103)
and a Key Program of NSF (Grant No. 10331020) in China.
\endthanks
\endtopmatter
\document

\heading{1. Introduction}\endheading

In a recent paper of the author and D. M. Davis [SD] originally motivated by the study of
homotopy exponents of the special unitary group $\t{SU}(n)$, the following sophisticated
theorem was established.

\proclaim{Theorem 1.0 {\rm (Sun and Davis)}} Let $p$ be a prime,
and let $\al,n\in\N=\{0,1,\ldots\}$
and $r\in\Z$. Then, for any $f(x)\in\Z[x]$, we have
 $$\aligned&\ord_p\(\sum_{k\eq r\,(\mo\ p^\al)}(-1)^k\bi nkf\l(\f{k-r}{p^\al}\r)\)
 \\\gs&\ord_p\l(\l\lfloor\f n{p^{\al-1}}\r\rfloor!\r)-\deg f
 +\tau_p(\{r\}_{p^{\al-1}},\{n-r\}_{p^{\al-1}}),
 \endaligned\tag1.0$$
 where $\ord_p(a)=\sup\{m\in\N:\, p^m\mid a\}$ is the $p$-adic order of $a\in\Z$,
 $\{a\}_{p^{\al-1}}$ stands for the least nonnegative residue of $a$ modulo $p^{\al-1}$
 $($and this is regarded as $0$ if $\al=0)$,
 and for $a,b\in\N$ we use $\tau_p(a,b)$ to denote the number of carries
 occurring in the addition of $a$ and $b$ in base $p$.
 \endproclaim

 Let $p$ be a prime. By a well-known fact in number theory (cf. [IR, p.\,26]),
$$\ord_p(n!)=\sum_{i=1}^{\infty}\l\lfloor \f n{p^i}\r\rfloor
\qquad \ \t{for every}\ n=0,1,2,\ldots.$$
A useful theorem of E. Kummer asserts that if $a,b\in\N$ then
$$\ord_p\bi{a+b}a=\sum_{i=1}^{\infty}\(\l\lfloor\f{a+b}{p^i}\r\rfloor
-\l\lfloor\f a{p^i}\r\rfloor-\l\lfloor\f b{p^i}\r\rfloor\)=\tau_p(a,b).$$

 In this paper we will apply Theorem 1.0 to
 deduce three theorems on combinatorial congruences
 or Stirling numbers of the second kind.

For $l,m\in\N$ with $l+m>0$,
the Stirling number $S(l,m)$ of the second kind
denotes the number of ways to partition a set of cardinality $l$ into
$m$ nonempty subsets; in addition, we define $S(0,0)$ to be $1$. It is well known that
$$x^l=\sum_{j=0}^{l}S(l,j)(x)_j\quad \ \t{for}\ l=0,1,2,\ldots,$$
where $(x)_j=\prod_{0\ls i<j}(x-i)$ and an empty product has the value $1$ (thus $(x)_0=1$).

Here is our first theorem.

\proclaim{Theorem 1.1} Let $p$ be any prime.
Let $a\in\Z$,  $l,l',m\in\Z^+=\{1,2,\ldots\}$, $l'\gs l>m/p$
and
$$l'\eq l\ \l(\mo\ (p-1)p^{\lfloor\log_p m\rfloor-\da_p(a,m)}\r),\tag1.1$$
where
$$\da_p(a,m)=\cases1&\t{if}\ a\in p\Z\ \t{and}\ \log_pm\in\Z^+,
\\0&\t{otherwise}.\endcases\tag1.2$$
Then we have
$$\sum_{j=0}^{l'}\bi {l'}jS(j,m)a^{l'-j}
\eq\sum_{j=0}^{l}\bi {l}jS(j,m)a^{l-j}\ (\mo\ p).\tag1.3$$
\endproclaim

\proclaim{Corollary 1.1} Let $p$ be a prime, and let $a\in\Z$ and $m\in\Z^+$.
Then, for $k=m+(p-1)p^{\lfloor\log_p m\rfloor-\da_p(a,m)}q$ with $q\in\N$, we have
$$\sum_{j=0}^{k}\bi kjS(j,m)a^{k-j}\eq 1\ (\mo\ p).$$
\endproclaim
\Proof. Just apply Theorem 1.1
with $l=m$ and $l'=k$. \qed

\Remark\ 1.1. Note that if  $p$ is a prime and $m$ is a positive integer then
$m-(p-1)p^{\lfloor\log_pm\rfloor}<m/p$.

\medskip

The following result was first obtained by
L. Carlitz [C] in 1955. (See also A. Nijenhuis and H. S. Wilf [NW],
and Y. H. H. Kwong [K].)

\proclaim{Corollary 1.2}  Let $p$ be any prime. Suppose that
$\al,m\in\N$, $m\gs p$ and $p^\al<m\ls p^{\al+1}$.
Then $p^\al(p-1)$ is a period of the sequence $\{S(l,m)\}_{l\gs m}$ modulo $p$.
\endproclaim
\Proof. It suffices to apply Theorem 1.1 with $a=0$. \qed

\medskip

The sum $\sum_{k\eq r\,(\mo\ m)}\bi nk$ with $m\in\Z^+$, $n\in\N$ and $r\in\Z$
has been investigated intensively, see [S] for some historical background
and related congruences.
In 1899 J.W.L. Glaisher (cf. [D, p.\,271] and [ST]) proved that
$$\sum_{j\eq r\,(\mo\ p-1)}\bi{l'}j\eq\sum_{j\eq r\,(\mo\ p-1)}\bi lj\ \ (\mo\ p)$$
whenever $p$ is a prime, $r\in\Z$, $l',l\in\Z^+$ and $l'\eq l\ (\mo\ p-1)$.
Clearly Glaisher's congruence is our following result in the case $m=1$.

\proclaim{Corollary 1.3} Let $p$ be a prime, $m\in\Z^+$ and $r\in\Z$.
For any $l',l\in\Z^+$ with $l'\gs l>m/p$ and
$$l'\eq l\ \l(\mo\ (p-1)p^{\lfloor\log_p m\rfloor}\r),$$
we have
$$\sum_{j\eq r\,(\mo\ p-1)}\bi{l'}jS(j,m)
\eq\sum_{j\eq r\,(\mo\ p-1)}\bi ljS(j,m)\ \ (\mo\ p).\tag1.4$$
\endproclaim

Our second theorem is slightly stronger than Conjecture 1.3 of the author and Davis [SD]
which was proved in [SD] when $p=2$ and $r=0$.

\proclaim{Theorem 1.2} Let $p$ be a prime, and let $\al,l,n\in\N$ and $r\in\Z$.
Set $r_*=\{r\}_{p^{\al}}$, $n_*=r_*+\{n-r\}_{p^{\al}}$ and
$$m=\f{n-n_*}{p^\al}=\l\lfloor \f r{p^{\al}}\r\rfloor
+\l\lfloor \f{n-r}{p^{\al}}\r\rfloor.\tag1.5$$
Suppose that $l\gs m>0$ and
$$l\eq m\ \l(\mo\ (p-1)p^{\lfloor\log_pm\rfloor
-\da_p(\lfloor r/p^\al\rfloor,\,m)}\r),\tag1.6$$
where the notation $\da_p(a,m)$ is given by $(1.2)$.
Then we have
$$\f1{\lfloor n/p^{\al}\rfloor!\bi{n_*}{r_*}}
\sum_{k\eq r\,(\mo\ p^{\al})}(-1)^k\bi nk\l(\f{k-r}{p^{\al}}\r)^{l}
\eq(-1)^{l+r_*}\ \ (\mo\ p).\tag1.7$$
\endproclaim

\Remark\ 1.2. Theorem 1.2 implies that the inequality in Theorem 5.1 of [DS]
is sharp for infinitely many values of $l$ provided that $n\gs 2p^\al-1$.
\medskip

 Our third theorem confirms Conjecture 1.1 of [SD].

 \proclaim{Theorem 1.3} Let $p$ be any prime, and let $\al\in\Z^+$, $l,n\in\N$
 and $r\in\Z$. Then
 $$\aligned &\f1{\lfloor n/p^{\al-1}\rfloor!}
\sum_{k\eq r\,(\mo\ p^{\al})}(-1)^{k}\bi{pn}{pk}\l(\f{k-r}{p^{\al-1}}\r)^l
\\\eq&\f1{\lfloor n/p^{\al-1}\rfloor!}
\sum_{k\eq r\,(\mo\ p^{\al})}(-1)^k\bi nk\l(\f{k-r}{p^{\al-1}}\r)^l\ (\mo\ p^{a_p}),
\endaligned\tag1.8$$
where
$$a_p=\cases 1&\t{if}\ p=2,\\2&\t{if}\ p=3,\\3&\t{if}\ p>3.
\endcases\tag1.9$$
\endproclaim

\Remark\ 1.3. Let $p$ be a prime, $\al,n\in\N$ and $r\in\Z$.
When $p^{\al}>n$ and $l=0\ls r\ls n$, (1.8) reduces to Ljunggren's
congruence $\bi{pn}{pr}\eq\bi nr\ (\mo\ p^{a_p})$
(cf. [G]) which is an extension of the Wolstenholme congruence
$\bi{2p}p\eq 2\ (\mo\ p^{a_p})$ (i.e., $\bi{2p-1}{p-1}\eq1\ (\mo\ p^{a_p})$).
Note also that (1.8) holds for every $l\in\N$ if and only if
we have
$$\bm{pn}{pr}{f,\,p^{\al+1}}\eq \bm nr{f,\,p^\al}\ \ (\mo\ p^{a_p})\tag1.10$$
for all $f(x)\in\Z[x]$, where
$$\bm nr{f,\,p^\al}=\f{p^{\deg f}}{\lfloor n/p^{\al-1}\rfloor!}
\sum_{k\eq r\,(\mo\ p^\al)}(-1)^k\bi nkf\l(\f{k-r}{p^\al}\r)\in\Z_p.\tag1.11$$
(As usual, $\Z_p$ denotes the ring of $p$-adic integers.)
\smallskip

Concerning the right-hand side of the congruence (1.8), a Lucas-type
congruence modulo $p$ was established in [SD] for $\al>1$ (and in [SW] for $\al=1$).
See also [SW] for some other congruences of Lucas' type related to combinatorial sums
involving binomial coefficients.

In the next section we are going to prove Theorem 1.1 and Corollary 1.3.
On the basis of Theorem 1.1 we will deduce Theorem 1.2 in Section 3.
Section 4 is devoted to our proof of Theorem 1.3.

\heading{2. Proofs of Theorem 1.1 and Corollary 1.3}\endheading

\noindent{\it Proof of Theorem 1.1}.
By a well-known property of Stirling numbers of the second kind (cf. [LW, pp.\,125--126]),
$$\align\sum_{j=0}^l\bi lja^{l-j}S(j,m)
=&\sum_{j=0}^l\bi lja^{l-j}\f1{m!}\sum_{k=0}^m(-1)^{m-k}\bi mkk^j
\\=&\f1{m!}\sum_{k=0}^m(-1)^{m-k}\bi mk(a+k)^l
=(-1)^m\sum_{r=0}^{p-1}S_r(l),
\endalign$$
where
$$S_r(l)=\f1{m!}\sum_{k\eq r\,(\mo\ p)}(-1)^{k}\bi mk(a+k)^l.\tag2.1$$

 Let $r\in\{0,\ldots,p-1\}$. Observe that
$$\align S_r(l)=&\f1{m!}\sum_{k\eq r\,(\mo\ p)}(-1)^{k}\bi mk
\sum_{j=0}^l\bi ljp^{j}\l(\f{k-r}{p}\r)^j(a+r)^{l-j}
\\=&\sum_{j=0}^l\bi lj (a+r)^{l-j}\f{p^{j}}{m!}
\sum_{k\eq r\,(\mo\ p)}(-1)^{k}\bi mk\l(\f{k-r}{p}\r)^j.
\endalign$$
By Theorem 1.0,  for any $j\in\N$ we have
 $$\sigma_r(j):=\f{p^j}{m!}\sum_{k\eq r\,(\mo\ p)}(-1)^{k}\bi mk\l(\f{k-r}p\r)_j\in\Z_p\tag2.2$$
and
 $$\f{p^j}{m!}\sum_{k\eq r\,(\mo\ p)}
 (-1)^k\bi mk\(\l(\f{k-r}p\r)^j-\l(\f{k-r}p\r)_j\)\eq0\ \ (\mo\ p)$$
since the degree of $f_j(x)=x^j-(x)_j\in\Z[x]$ is smaller than $j$.
Therefore,
$$S_r(l)\eq\sum_{j=0}^l\bi lj (a+r)^{l-j}\sigma_r(j)\ \ (\mo\ p).\tag2.3$$

In view of the above, it suffices to show that
$$\bi {l'}j (a+r)^{l'-j}\sigma_r(j)\eq\bi lj (a+r)^{l-j}\sigma_r(j)\ \ (\mo\ p)\tag2.4$$
for every $j=0,1,\ldots$.

Below we assume $j\in\N$ and $\sigma_r(j)\not=0$. Then
 $((k-r)/p)_j\not=0$ for some $0\ls k\ls m$ with $k\eq r\,(\mo\ p)$,
hence $m-r\gs k-r\gs pj$ and $j\ls m/p<l\ls l'$.
If $p\mid a+r$, then
$$(a+r)^{l'-j}\eq0\eq (a+r)^{l-j}\ \ (\mo\ p).$$
When $p\nmid a+r$, as $l'\eq l\ (\mo\ p-1)$ we have
$$(a+r)^{l'-j}\eq (a+r)^{l-j}\ \ (\mo\ p)$$
by Fermat's little theorem.
So it remains to show $\bi {l'}j\eq\bi lj\ (\mo\ p)$
in the case $p\nmid a+r$.

Let $\al=\lfloor\log_p m\rfloor$.
Then $p^{\al}\ls m<p^{\al+1}$ and $\da=\da_p(a,m)\ls\al$.
Write $l=p^{\al-\da}q_0+l_0$ with $q_0\in\N$
and $0\ls l_0<p^{\al-\da}$. For some $q\in\N$ we have
$l'=l+(p-1)p^{\al-\da}q=p^{\al-\da}((p-1)q+q_0)+l_0$.
Recall that $j\ls(m-r)/p<p^\al$.
Suppose $a+r\not\eq0\ (\mo\ p)$.
If $\da=1$, then $j<m/p=p^{\al-1}$ because
$m=p^\al$ and $r\not=\{-a\}_p=0$. Thus
$j<p^{\al-\da}$. With help of the Chu-Vandermonde convolution identity (cf. [GKP, (5.27)]),
$$\align \bi lj-\bi{l_0}j=&\sum_{0<i\ls j}\bi{p^{\al-\da}q_0}{i}\bi{l_0}{j-i}
\\=&\sum_{0<i\ls j}\f{p^{\al-\da}}iq_0\bi{p^{\al-\da}q_0-1}{i-1}\bi{l_0}{j-i}\eq0\ (\mo\ p).
\endalign$$
Similarly, $\bi {l'}j\eq\bi{l_0}j\ (\mo\ p)$ as desired.
We are done. \qed

\medskip
\noindent{\it Proof of Corollary 1.3}. Let $g$ be a primitive root modulo $p$.
For any integer $h$, if $p-1\mid h$ then $\sum_{a=1}^{p-1}a^h\eq p-1\eq-1\ (\mo\ p)$
by Fermat's little theorem; if $p-1\nmid h$ then $g^h\not\eq1\ (\mo\ p)$
and hence $\sum_{a=1}^{p-1}a^h\eq0\ (\mo\ p)$ since
$$(g^h-1)\sum_{a=1}^{p-1}a^h=\sum_{a=1}^{p-1}(ag)^h-\sum_{a=1}^{p-1}a^h\eq0\ (\mo\ p).$$

In view of the above,
$$\align&\sum_{a=1}^{p-1}a^{r-l}\sum_{j=0}^l\bi ljS(j,m)a^{l-j}
\\=&\sum_{j=0}^l\bi ljS(j,m)\sum_{a=1}^{p-1}a^{r-j}
\\\eq&-\sum_{j\eq r\,(\mo\ p-1)}\bi ljS(j,m)\ (\mo\ p).
\endalign$$
Similarly,
$$\sum_{a=1}^{p-1}a^{r-l'}\sum_{j=0}^{l'}\bi {l'}jS(j,m)a^{l'-j}
\eq-\sum_{j\eq r\,(\mo\ p-1)}\bi {l'}jS(j,m)\ \ (\mo\ p).$$
Since $l'\eq l\ (\mo\ p-1)$, $a^{r-l'}\eq a^{r-l}\ (\mo\ p)$
for all $a=1,\ldots,p-1$. Thus, applying Theorem 1.1
we immediately obtain (1.4) from the above. \qed

\heading{3. Proof of Theorem 1.2}\endheading

At first we make some useful observations. Clearly
$$m\ls\f{n-r_*}{p^\al}=\f{n-r}{p^\al}+\l\lfloor\f r{p^{\al}}\r\rfloor
<1+\l\lfloor\f{n-r}{p^\al}\r\rfloor
+\l\lfloor\f r{p^\al}\r\rfloor=m+1.$$
Since
$$\tau_p(\{r\}_{p^\al},\{n-r\}_{p^\al})-\tau_p(\{r\}_{p^{\al-1}},\{n-r\}_{p^{\al-1}})
=\cases1&\t{if}\ n_*\gs p^\al,
\\0&\t{otherwise},\endcases$$
we also have
$$\ord_p\bi{n_*}{r_*}-\tau_p(\{r\}_{p^{\al-1}},\{n-r\}_{p^{\al-1}})
=\l\lfloor \f{n_*}{p^\al}\r\rfloor=\l\lfloor\f n{p^\al}\r\rfloor-m.\tag3.1$$

Let
$a=-\lfloor r/p^{\al}\rfloor$, $k\in\{0,\ldots,n\}$ and $k\eq r\ (\mo\ p^\al)$. Then
$$\align \l(\f{k-r}{p^\al}\r)^l=&\l(\f{k-r_*}{p^\al}+a\r)^l
=\sum_{j=0}^l\bi lja^{l-j}\l(\f{k-r_*}{p^\al}\r)^j
\\=&\sum_{j=0}^l\bi lj a^{l-j}\sum_{i=0}^jS(j,i)\l(\f{k-r_*}{p^\al}\r)_i
\\=&\sum_{j=0}^l\bi lj a^{l-j}\sum_{i=0}^mS(j,i)\l(\f{k-r_*}{p^\al}\r)_i,
\endalign$$
because for $i\gs m+1$ we have
$i>(n-r_*)/p^\al\gs(k-r_*)/p^\al$ and hence $((k-r_*)/p^\al)_i=0$.

Observe that
$$\ord_p\(\l\lfloor\f n{p^{\al-1}}\r\rfloor!\)
=\l\lfloor\f n{p^\al}\r\rfloor+\sum_{s>\al}\l\lfloor\f n{p^s}\r\rfloor
=\l\lfloor\f n{p^\al}\r\rfloor+\ord_p\(\l\lfloor\f n{p^{\al}}\r\rfloor!\).$$
If $i\in\{0,\ldots,m-1\}$, then by Theorem 1.0 and (3.1) we have
$$\align&\ord_p\(\sum_{k\eq r_*\,(\mo\ p^\al)}(-1)^k\bi nk\l(\f{k-r_*}{p^\al}\r)_i\)
\\\gs&\ord_p\(\l\lfloor\f n{p^{\al-1}}\r\rfloor!\)-i
+\tau_p(\{r_*\}_{p^{\al-1}},\{n-r_*\}_{p^{\al-1}})
\\>&\ord_p\(\l\lfloor\f n{p^{\al}}\r\rfloor!\)+\l\lfloor\f n{p^\al}\r\rfloor-m
+\tau_p(\{r\}_{p^{\al-1}},\{n-r\}_{p^{\al-1}})
\\=&\ord_p\(\l\lfloor\f n{p^\al}\r\rfloor!\bi{n_*}{r_*}\).
\endalign$$
Therefore,
$$\align&\f1{\lfloor n/p^\al\rfloor!\bi{n_*}{r_*}}
\sum_{k\eq r\,(\mo\ p^\al)}(-1)^k\bi nk\l(\f{k-r}{p^\al}\r)^l
\\\eq&\sum_{j=0}^l\bi lj\f{S(j,m)a^{l-j}}{\lfloor n/p^\al\rfloor!\bi{n_*}{r_*}}
\sum_{k\eq r\,(\mo\ p^\al)}(-1)^k\bi nk\l(\f{k-r_*}{p^\al}\r)_m\ (\mo\ p).
\endalign$$
In light of Corollary 1.1, it remains to show that
$S\eq(-1)^{l+r_*}\ (\mo\ p)$, where
$$S=\f1{\lfloor n/p^\al\rfloor!\bi{n_*}{r_*}}
\sum_{k\eq r_*\,(\mo\ p^\al)}(-1)^k\bi nk\l(\f{k-r_*}{p^\al}\r)_m.\tag3.2$$

If $k\in\{0,\ldots,n\}$, $k\eq r_*\ (\mo\ p^\al)$ and $((k-r_*)/p^\al)_m\not=0$, then
$$m\ls\f{k-r_*}{p^\al}\ls\f{n-r_*}{p^\al}<m+1$$
and hence $k=mp^\al+r_*$. So
$$S=\f{(-1)^{mp^\al+r_*}}{\lfloor n/p^\al\rfloor!\bi{n_*}{r_*}}\bi n{mp^\al+r_*}(m)_m
=\f{(-1)^{mp^\al+r_*}m!}{(m+\lfloor n_*/p^\al\rfloor)!}
\times\f{\bi {mp^\al+n_*}{mp^\al+r_*}}{\bi{n_*}{r_*}}.$$
Clearly
$$\align\f{\bi {mp^\al+n_*}{mp^{\al}+r_*}}{\bi{n_*}{r_*}}
=&\f{(mp^\al+n_*)!/(mp^{\al}+r_*)!}{n_*!/r_*!}
\\=&\prod_{0<i\ls n_*}\l(1+m\f{p^\al}i\r)\bigg/\prod_{0<j\ls r_*}\l(1+m\f{p^\al}j\r).
\endalign$$
Thus, if $n_*<p^\al$ then
$$\f{\bi {mp^\al+n_*}{mp^{\al}+r_*}}{\bi{n_*}{r_*}}\eq1\ (\mo\ p);$$
if $n_*\gs p^\al$ then $\lfloor n_*/p^\al\rfloor=1$ and
$$\f{\bi {mp^\al+n_*}{mp^{\al}+r_*}}{(m+1)\bi{n_*}{r_*}}
=\prod\Sb 0<i\ls n_*\\i\not=p^\al\endSb\l(1+m\f{p^\al}i\r)
\bigg/\prod_{0<j\ls r_*}\l(1+m\f{p^\al}j\r)
\eq1\ (\mo\ p).$$
Therefore,
$$S\eq(-1)^{mp^\al+r_*}\eq(-1)^{m+r_*}\eq(-1)^{l+r_*}\ \ (\mo\ p).$$
This concludes the proof of Theorem 1.2.

\heading{4. Proof of Theorem 1.3}\endheading

For $i,k\in\N$ let $\da_{i,k}$ be the Kronecker symbol
which takes $1$ or $0$ according to whether $i=k$ or not.
Since
$$\da_{i,k}=\bi ki\sum_{j\gs i}(-1)^{j-i}\bi{k-i}{j-i}=\sum_{j\gs i}(-1)^{j-i}\bi{k}j\bi ji,$$
we have
$$\align&(-1)^{(p-1)r}\sum_{k\eq r\,(\mo\ p^{\al})}
(-1)^{k}\bi {pn}{pk}\l(\f{k-r}{p^{\al-1}}\r)^l
\\=&\sum_{k\eq r\,(\mo\ p^{\al})}(-1)^{pk}\bi {pn}{pk}\l(\f{k-r}{p^{\al-1}}\r)^l
\\=&\sum_{k=0}^n(-1)^{pk}\bi{pn}{pk}\sum_{i\eq r\,(\mo\ p^{\al})}
\l(\f{i-r}{p^{\al-1}}\r)^l
\da_{i,k}
\\=&\sum_{k=0}^n(-1)^{pk}\bi{pn}{pk}\sum_{i\eq r\,(\mo\ p^{\al})}
\l(\f{i-r}{p^{\al-1}}\r)^l
\sum_{j\gs i}(-1)^{j-i}\bi{k}j\bi ji
\\=&\sum_{j=0}^n(-1)^jC_{n,j}
\sum_{i\eq r\,(\mo\ p^{\al})}(-1)^i\bi ji\l(\f{i-r}{p^{\al-1}}\r)^l,
\endalign$$
where
$$C_{n,j}=\sum_{k=0}^n(-1)^{pk}\bi{pn}{pk}\bi{k}j
=\sum_{k\eq 0\,(\mo\ p)}(-1)^k\bi {pn}k\bi{k/p}j.$$

 As $C_{n,n}=(-1)^{pn}$, by the above
$$\align &(-1)^{(p-1)r}\sum_{k\eq r\,(\mo\ p^{\al})}(-1)^{k}
\bi{pn}{pk}\l(\f{k-r}{p^{\al-1}}\r)^l
\\&-(-1)^{(p-1)n}\sum_{i\eq r\,(\mo\ p^\al)}(-1)^i\bi ni\l(\f{i-r}{p^{\al-1}}\r)^l
\\=&\sum_{0\ls j<n}(-1)^jC_{n,j}
\sum_{i\eq r\,(\mo\ p^{\al})}(-1)^i\bi ji\l(\f{i-r}{p^{\al-1}}\r)^l.
\endalign$$
Note that $(-1)^{(p-1)n}\eq(-1)^{(p-1)r}\ (\mo\ p^{a_p})$.
In view of Theorem 1.0,
$$\ord_p\(\sum_{i\eq r\,(\mo\ p^\al)}(-1)^i\bi ji\l(\f{i-r}{p^{\al-1}}\r)^l\)
\gs\ord_p\(\l\lfloor\f j{p^{\al-1}}\r\rfloor!\)
=\sum_{s=\al}^{\infty}\l\lfloor\f j{p^s}\r\rfloor.$$
So it suffices to show that
$$\ord_p(C_{n,j})\gs a_p+\ord_p\(\l\lfloor\f n{p^{\al-1}}\r\rfloor!\)
-\sum_{s=\al}^{\infty}\l\lfloor\f j{p^s}\r\rfloor
=a_p+\sum_{s=\al}^{\infty}\(\l\lfloor\f n{p^s}\r\rfloor-\l\lfloor\f j{p^s}\r\rfloor\)$$
for any $j\in\N$ with $j<n$.

Fix a nonnegative integer $j<n$.
In light of Theorem 1.0,
$$\ord_p(j!C_{n,j})\gs\ord_p\l(\l\lfloor\f{pn}{p^{1-1}}\r\rfloor!\r)-j
=\sum_{s=0}^{\infty}\l\lfloor\f n{p^s}\r\rfloor-j.$$
By Lemma 3.2 of [SD] and its proof, $C_{n,j}$ is congruent to
$$\sum_{k=0}^n(-1)^{k}\bi{n}{k}\bi kj=(-1)^j\bi nj\sum_{k\gs j}(-1)^{k-j}\bi{n-j}{k-j}
=0$$
modulo $p^{2\ord_p(n)+a_p}$.
In the case $p>3$, by Jacobsthal's result (cf. [G]), if $k\in\{1,\ldots,n\}$, then
$$\bi{pn}{pk}\bigg/\bi nk=1+p^3nk(n-k)q_k$$
for some $q_k\in\Z_p$,
and hence
$$\bi{pn}{pk}-\bi nk=\bi nkp^3nk(n-k)q_k=p^3n^2(n-1)\bi{n-2}{k-1}q_k.$$
So we also have $\ord_p(C_{n,j})\gs \ord_p(n-1)+3$ when $p>3$.
These facts will be used in the following discussion.

{\it Case} 1. $n-j\gs  a_p$.

In this case,
$$\align \ord_p(C_{n,j})\gs&\sum_{s=0}^{\infty}\l\lfloor\f n{p^s}\r\rfloor-j-\ord_p(j!)
=\sum_{s=0}^{\infty}\(\l\lfloor\f n{p^s}\r\rfloor-\l\lfloor\f j{p^s}\r\rfloor\)
\\=&n-j+\sum_{s=1}^{\infty}\(\l\lfloor\f n{p^s}\r\rfloor-\l\lfloor\f j{p^s}\r\rfloor\)
\\\gs&a_p+\sum_{s=\al}^{\infty}\(\l\lfloor\f n{p^s}\r\rfloor-\l\lfloor\f j{p^s}\r\rfloor\).
\endalign$$

{\it Case} 2. $0<n-j<a_p\ls3$, and $p\mid n$ or $j\not=n-2$.

If $\beta=\ord_p(n)>0$, then $n-j<a_p<p\ls p^\beta$ and hence
$$\l\lfloor\f n{p^{\beta+1}}\r\rfloor=\l\lfloor\f{n/p^\beta}{p}\r\rfloor
=\l\lfloor\f{n/p^\beta-1}p\r\rfloor=\l\lfloor\f{\lfloor j/p^\beta\rfloor}p\r\rfloor
=\l\lfloor\f j{p^{\beta+1}}\r\rfloor,$$
therefore
$$\align&\sum_{s=\al}^{\infty}\(\l\lfloor\f n{p^s}\r\rfloor-\l\lfloor\f j{p^s}\r\rfloor\)
=\sum_{\al\ls s\ls \beta}\(\f n{p^s}-\l\lfloor\f j{p^s}\r\rfloor\)
\\\ls&\sum_{\al\ls s\ls \beta}1<2\beta=2\ord_p(n)\ls \ord_p(C_{n,j})-a_p.
\endalign$$
When $\beta=\ord_p(n)=0$ (i.e., $p\nmid n$) and $j=n-1$,
we have
$$\sum_{s=\al}^{\infty}\(\l\lfloor\f n{p^s}\r\rfloor-\l\lfloor\f j{p^s}\r\rfloor\)
=0=2\ord_p(n)\ls \ord_p(C_{n,j})-a_p.$$

{\it Case} 3. $n-j=2<a_p$ and $p\nmid n$.

In this case, $a_p=3<p$ and
$$\align \ord_p(C_{n,j})-a_p\gs&\ord_p(n-1)
\gs\sum_{\al\ls s\ls\ord_p(n-1)}\(\f{n-1}{p^s}-\l\lfloor\f{n-2}{p^s}\r\rfloor\)
\\=&\sum_{s=\al}^{\infty}
\(\l\lfloor\f{n-1}{p^s}\r\rfloor-\l\lfloor\f{n-2}{p^s}\r\rfloor\)
=\sum_{s=\al}^{\infty}\(\l\lfloor\f{n}{p^s}\r\rfloor-\l\lfloor\f{j}{p^s}\r\rfloor\).
\endalign$$

Combining the above we have completed the proof of Theorem 1.3.

\bigskip

\leftline{Department of Mathematics} \leftline{Nanjing University}
\leftline{Nanjing 210093} \leftline{People's Republic of China}
\leftline {zwsun\@nju.edu.cn} \leftline {\tt
http://math.nju.edu.cn/${}^\sim$zwsun}

\enddocument